\documentclass[12pt]{article}
\usepackage{amssymb}
\usepackage[leqno]{amsmath}
\usepackage{amsfonts}
\usepackage{amsopn}
\usepackage{amstext}
\usepackage{amsthm}
\newcommand{\subs}{\subseteq}
\def\Proof{{\it Proof. }}

\def\eqdef{\mbox{\bf\ :=\ }}

\newcommand{\rfs}[1]{{\rm #1}}
\newcommand{\fnn}[3]{#1:#2 \rightarrow #3}

\newcommand{\setn}[2]{\{#1\: :\:#2\}}

\newcommand{\setd}[2]{\{#1,\dots,#2\}}
\newcommand{\seqn}[2]{\langle#1\: :\:#2\rangle }
\newcommand{\sngltn}[1]{ \{#1\} }
\newcommand{\dbltn}[2]{\{ #1, #2 \}}
\newcommand{\pair}[2]{\langle #1 , \:#2 \rangle }
\newcommand{\triple}[3]{\langle\kern1pt#1 , \:#2 , 
      \:#3 \kern1pt\rangle }
\newcommand{\forthuple}[4]{\langle\kern1pt#1 , \:#2  , 
      \:#3  , \:#4\kern1pt\rangle }
\newcommand{\qdrpl}[4]{\langle\kern1pt#1, #2, #3, 
   #4\kern1pt\rangle}
\newcommand{\fifthtpl}[5]{\langle\kern1pt#1, #2, #3, 
   #4, #5\kern1pt\rangle}
\newcommand{\sixtpl}[6]{\langle\kern1pt#1, #2, #3, 
   #4, #5, #6 \kern1pt\rangle}

\newtheorem{definition}{{Definition}}[section]
\newtheorem{theorem}[definition]{{Theorem}}

\newtheorem{proposition}[definition]{{Proposition}}
\newtheorem{lemma}[definition]{\noindent {Lemma}}
\newtheorem{question}[definition]{\noindent {Question}}

\newtheorem{fact}[definition]{\noindent {Fact}}

\newtheorem{remark}[definition]{\noindent{Remark}}

\newcommand{\meet}{\wedge}
\newcommand{\join}{\vee}

\providecommand{\cal}{\mathcal}

\def\natN{{\mathbb{N}\kern1pt}}
\def\natQ{\hbox{$\mathbb{Q}\kern1pt$}}

\newenvironment{pf}{\begin{proof}}{\end{proof}}

\def\calA{\hbox{$\mathcal{A}\kern1pt$}}
\def\calB{\mathcal{B}}
\def\calG{\hbox{\kern2pt$\mathcal{G}$}}
\def\calP{\hbox{\kern2pt$\mathcal{P}$}}

\def\calG{{\mathcal{G}}}
\def\calS{{\mathcal{S}}}

\def\F{F}
\def\L{L}
\def\Fs{\rfs{FS}}
\def\Is{\rfs{IS}}
\def\Pr{\Pi}
\def\classw{{\calB}}

\def\gwf{{\calG}_{wf}}
\def\gwqo{{\calG}_{wqo}}
\def\gbqo{{\calG}_{bqo}}

\def\BQO{{\cal BQO}}
\def\WQO{{\cal WQO}}

\newcommand{\tleq}{\triangleleft}

\newcommand{\ignore}[1]{} 

\begin{document}
%
%
\thispagestyle{empty}
\phantom{abk}
\vspace{-5mm}
\begin{center}
{\Large\bf Poset algebras over well quasi-ordered posets\vspace{3mm}}

{\bf Uri Abraham}

{\small Department of Mathematics,\\
Ben Gurion University,
Beer-Sheva, Israel\vspace{1mm}
}

{\bf Robert Bonnet}%
\footnote{This work is supported by the Center for Advanced Studies
in Mathematics (Ben Gurion University).}

{\small Laboratoire de Math\'ematiques,\\
Universit\'e de Savoie, Le Bourget-du-Lac, France \vspace{1mm}
}

{\bf Wies{\l}aw Kubi\'s}%
\footnote{This work is supported by 
the Israel Science Foundation (Post-Doctoral positions
at Ben Gurion University 2000--2002),
the Fields Institute (Toronto 2002--2004),
and by the Nato Science Fellowship
(University Paris~VII, CNRS-UMR 7056, 2004).}

{\small 
Instytut Matematyki,\\
Akademia \'Swi\c{e}tokrzyska,
Kielce, Poland \vspace{1mm}
}

{\bf (20 February 2007)}
\end{center}

\begin{abstract}
A new class of partial order-types, class $\gbqo^+$ is defined and investigated here.
A poset $P$ is in the class $\gbqo^+ $ iff the poset algebra $\F(P)$ is generated 
by a better quasi-order $G$ that is included in $\L(P)$.

The free Boolean algebra $\F(P)$ and its free distrivutive lattice $\L(P)$ 
were defined in~\cite{ABKR}.
The free Boolean algebra $\F(P)$ contains the partial order $P$ and is generated by it:
$\F(P)$ has the following universal property.
If $B$ is any Boolean algebra and $f$ is any order-preserving map from $P$ 
into a Boolean algebra $B$, then $f$ can be extended  to an homomorphism $\hat{f}$
of $F(P)$ into $B$.
We also define $\L(P)$ as the sublattice of $F(P)$ generated by~$P$.

We prove that if $P$ is any well quasi-ordering, 
then $\L(P)$ is well founded, and is a countable union of well quasi-orderings.

We prove that the class $\gbqo^+$ is contained in the class of well quasi-ordered sets.
We prove that $\gbqo^+$ is preserved under homomorphic image, finite products, and lexicographic sum over better quasi-ordered index sets.
We prove also that every countable well quasi-ordered set is in~$\gbqo^+$.
We do not know, however if the class of well quasi-ordered sets is contained in~$\gbqo^+$.
Additional results concern homomorphic images of posets algebras.
\end{abstract}
\begin{footnotesize}
\newline
{\bf Mathematics Subject  Classification 2000 (MSC2000)}
\newline
{\bf Primary:}
03G05,  
06A06.  
{\bf Secondary:}
06E05,  
08A05,  
54G12.  
\medskip
\newline
{\bf E-mail:} abraham@math.bgu.ac.il (U. Abraham),
bonnet@in2p3.fr (R. Bonnet),
\\
wkubis@pu.kielce.pl (W. Kubi\'s)
\end{footnotesize}
\section{Introduction}
\label{S1}
Three important classes of partially ordered sets that interest us in this paper 
are the well-founded, well quasi-ordered, and better quasi-ordered classes. 
Intuitively, we tend to view these classes as good properties that a poset might posses. 
A poset not in a good class might nevertheless have a certain affinity with that class 
which makes it for that reason interesting. 
For example, $P^\ast$, the inverse ordering of $P$, may be well-founded, 
or $P$ may be a countable union of well-founded posets etc.  
Here we are interested in a different type of affinity of a poset $P$, 
namely that its poset algebra $F(P)$ is generated by a sublattice belonging to the good class.
The main class of posets introduced and studied here is the class ${\gbqo}$ 
of those posets that have this type of afinity with the better quasi-ordering.  
To understand this notion we must first recall the definition of the poset algebra 
defined in \cite{ABKR} and review some of the main results concerning 
posets algebras proved in that paper.
So we begin our introduction with some definitions and useful facts.

We recall some notions  concerning posets (partially ordered sets).
Let $\pair{P}{\leq^P}$ and $\pair{Q}{\leq^Q}$ be  posets. 
A map $f:P\to Q$ is {\em order-preserving} 
if $p \leq^P q$ implies that  $f(p)\leq^{Q} f(q)$. 
(When $\leq^{P}$ is understood from the context, we omit it and write instead $\leq$\,.)
A poset $P$ is {\it well founded\/} if there is no infinite descending chain
$p_0>p_1>p_2>\cdots$ in $P$.
Elements $p,q \in P$ are {\it incomparable} if neither $p\leq q$ nor $q\leq p$.
An {\it antichain} in a poset $P$ is a subset of $P$ consisting of pairwise 
incomparable elements. 
A poset $P$ is {\it narrow} if all antichains in $P$ are finite.
A poset $P$ is {\it well quasi-ordered} ({\it wqo\/}) if it is well founded and narrow.
A poset $P$ is {\it scattered} if it does not contain an isomorphic copy of the rational 
numbers~$\natQ$.

For Boolean algebras we use the notations of~\cite{kopp}.
Thus $+$, $\cdot$\,, $-$ and $\leq$ denote the join, meet,
complementation and partial ordering of a Boolean algebra~$B$.

A {\it homomorphism\/} $\fnn{g}{ B_1 }{ B_2 }$  of a Boolean algebra $B_1$ 
into a Boolean algebra $B_2$ is a function (not necessarily one-to-one) 
that respects the join, meet, and complementation operations.

Let $B$ be a Boolean algebra.
We say that $L\subseteq B$ is a {\it sublattice\/} if $L$ is closed under the meet and join 
operations of $B$.
The sublattice generated by a subset $X$ of $L$ is the minimal sublattice of $B$ that 
includes $X$. 
The Boolean subalgebra generated by $X$ is the minimal subalgebra of $B$ that contains $X$.

Following~\cite{ABKR}, we say that a Boolean algebra $B$ is {\it well-generated\/} 
if $B$ has a well-founded sublattice $G$ that generates $B$ 
(that is, $G$ is well-founded under the ordering of $B$ and the subalgebra generated 
by $G$ is $B$).

Similarly, we say that $B$ is {\it wqo-generated} ({\it better-generated\/\/}) 
if $B$ has a sublattice that generates $B$ and is a well quasi-order 
(better quasi-order) under the ordering of $B$. For the definition of a better quasi-ordering we refer to Section \ref{S2}.

Let $P$ be a poset.
The {\it poset algebra\/} $\F(P)$ is a Boolean algebra
(that turns out to be unique) that satisfies the following.
There exists an order-preserving injection $\fnn{x}{ P }{ \F(P) }$ such that:
\begin{enumerate}
\item The image  $x[P]$ of $x$ generates $\F(P)$ and

\item if $B$ is any Boolean algebra and $\fnn{f}{ P }{ B }$ any order-preserving map 
(not necessarily an injection) then there exists 
a Boolean homomorphism $\fnn{ \hat{f} }{ \F(P) }{ B }$ such that $\hat{f}(x_p) = f(p)$ 
for every $p \in P$, where $x_p$ is just a way of writing $x(p)$.
So $f =  \hat{f} \circ x$.
\end{enumerate}
We also define $\L(P)$ as the sublattice of $F(P)$ generated by 
$x[P] = \setn{ x_p }{ p \in P }$.

For a poset $P$,
we denote by $\Pr(P)$ the meet subsemilattice of $\F(P)$ generated by $x[P]$.
So a member of $\Pr(P)$ is of the form
$x_{\sigma} \eqdef \prod\setn{ x_p }{ p \in \sigma }$ where $\sigma$ is a finite subset of $P$;
and $\L(P)$ is the join subsemilattice of $\F(P)$ generated by $\Pr(P)$.
So $\Pr(P) \subseteq \L(P) \subseteq \F(P)$.

The {\it interval algebra\/} of a linearly ordered set $\pair{L}{\leq}$
is the subalgebra of $\calP(L)$ generated by the family $\setn{ [a,\rightarrow) }{a \in L}$ 
of left-closed rays of $L$, where $[a,\rightarrow)= \setn{x\in L}{a\leq x}$.
This algebra is denoted by $B(L)$.
The interval algebra of $L$ is isomorphic to the poset algebra of $L$ 
when $L$ has no minimum, and to the poset algebra of $L$ minus its minimum if $L$
has a minimum.
So the class of poset algebras contains the class of interval algebras.
\medskip

A Boolean algebra $B$ is {\it superatomic\/},
if every homomorphic image of~$B$ has an atom.
\vspace{-2mm}
\begin{itemize}
\item[($\star$)]
{\it Every well generated algebra is superatomic.}
\end{itemize}
\vspace{-2mm}
For the proof we refer to \cite[Proposition 2.7(b)]{BR1}.
\smallskip

The notion of poset algebra enables the following three definitions:
\begin{definition}
\label{dfn-1.1}
\begin{rm}
\begin{enumerate}
\item 
The class ${\gwf}$ contains all the posets $P$ such that $\F(P)$ is well-generated.
\item 
The class $\gwqo$ contains all the posets $P$ such that $\F(P)$ is 
well quasi-ordered generated.
\item 
The class ${\gbqo}$ contains all the posets $P$ such that $\F(P)$ is better generated.
\end{enumerate}
\end{rm}
\end{definition}

The class ${\gwf}$ was completely characterized in 
~\cite[Theorems~1.3]{ABKR}  as follows.
\medskip
\newline
{\bf Theorem A}
\begin{it}
Let $P$ be a poset. The following conditions are equivalent:

\rfs{(i)} $P\in \gwf$, that is $\F(P)$ is well generated.

\rfs{(ii)} $P$ is scattered and narrow.

\rfs{(iii)} $\F(P)$ is superatomic.
\end{it}
\medskip
\newline
In fact, M. Pouzet has proved (in an earlier work) that condition~(ii) 
is equivalent to condition~(iii).
Recall that by~($\star$), (ii) implies (iii). 
But there are superatomic Boolean algebras 
which are not well generated~\cite[Theorem~3.4]{BR1}.

If $P$ is a well quasi-ordered poset, then $P$ is certainly scattered and narrow, 
and hence $\F(P)$ is well-generated. 
We improve this result and prove in Theorem~\ref{thm-2.5} the following fact.
Suppose that $W$ is a well quasi-ordering.
Then the lattice $\L(W)$ generated by $x[W] = \setn{ x_p }{ p \in W }$
has the following properties:
\vspace{-1mm}
\begin{itemize}
\item[($\dagger$)] $\L(W)$ is a countable union of well quasi-orderings.
\vspace{-2mm}
\item[($\ddagger$)] $\L(W)$ is well founded.
\vspace{-1mm}
\end{itemize}
In particular, by ($\dagger$),
every antichain of~$\L(W)$ is countable.

It is impossible to improve this result and obtain that $\L(W)$ is well quasi-orderred, 
because the poset $R$ of Rado~\cite{rado} would be a counterexample. 
This poset is a countable well quasi-ordering that is not a better quasi-ordering and 
we shall observe that $\L(R)$ has an infinite antichain (see Remark~\ref{remark-2.6}).
On the other hand, we shall prove that $\F(R)$ is generated by a better quasi-ordered 
lattice contained in $\L(R)$ (see Theorem~\ref{thm-2.12} below).

We tend to believe that well quasi-ordering have greater affinity to 
better quasi-ordering than the definitions would allow us to think.  
We make the following two conjectures:
\begin{enumerate}
\item Any well quasi-ordering is a countable union of better quasi-orderings.
\item If $P$ is a well quasi-ordering, then $P\in \gbqo$ (or at least $P\in \gwqo$).
\end{enumerate}
Better quasi-orderings have, as the name suggests,
nicer properties than well quasi-ordered sets; for example: 
if $Q$ is a better quasi-ordered subset of a Boolean algebra $B$ then 
the sublattice of $B$ generated by $Q$ is again better quasi-ordered 
(Proposition~\ref{prop-2.3}(c))
---this is not true for well quasi-orderings (consider the Rado poset).
Thus, a Boolean algebra which is generated by  a better quasi-ordered subset is also 
generated by a better quasi-ordered sublattice and therefore it is well generated.

\begin{figure}
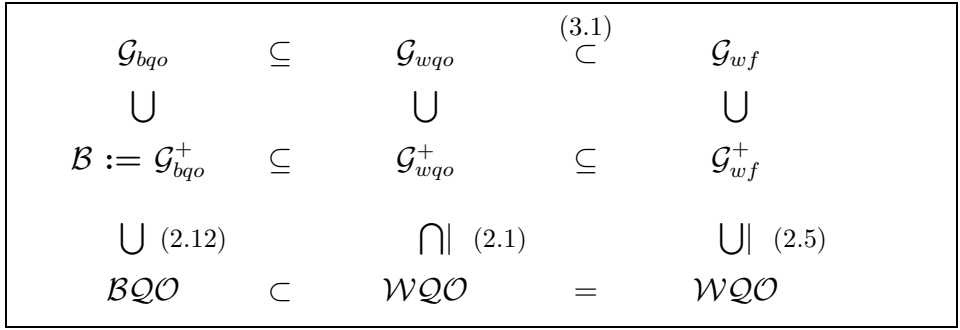
\label{figure-one}
\fbox{
\begin{tabular}{cccccccl}
\vspace{2mm} 
$\gbqo$ & $\subseteq$ &  $\gwqo$ &  
${\buildrel\mbox{\rm \footnotesize (\ref{thm-3.1})} \over \subset}$ & $\gwf$ & 
 \\
\vspace{2mm} 
$\bigcup$ &     & $\bigcup$  &    &  $\bigcup$   & & &  \\
$\calB \eqdef \gbqo^+$ & $\subseteq$ &  $\gwqo^+$ &  $\subseteq$ & $\gwf^+$  & \\ \\
\vspace{2mm} 
$\phantom{ \mbox{\rm\footnotesize\ (\ref{thm-2.1})} }
\bigcup
\mbox{\rm\footnotesize\ (\ref{thm-2.12}) }$ 
&     
& \phantom{ \mbox{\rm\footnotesize\ (\ref{thm-2.1})} } 
$\bigcap\!\vert$ \mbox{\rm\footnotesize\ (\ref{thm-2.1})} &    
                    &  \phantom{\mbox{\rm\footnotesize\ (\ref{thm-2.5})}}
                    $\bigcup\!\vert$ \mbox{\rm\footnotesize\ (\ref{thm-2.5})} &    &  & \\
\vspace{1mm} 
$\BQO$   &  
${\subset}$   
&  $\WQO$  &  $=$   & $\WQO$     &     &  & \\
\end{tabular}
}
\caption{Relations between the classes. Recall that $\gwf$  
=   scattered and narrow. }
\end{figure}

Concerning the second conjecture above, we prove in Theorem~\ref{thm-2.12}
that it holds for every countable well quasi-ordering. 
In fact, we prove a stronger result for these orderings: 
not only that the poset algebra of a countable well quasi-ordering is better generated, 
but there is a better quasi-ordered sublattice of $\L(P)$ that generates $\F(P)$.
This leads to the following definitions.
\begin{enumerate}
\item 
The class $\gwf^{+}$ contains all the posets $P$ such that $\F(P)$ is 
\emph{positively} well-generated, 
that is a sublattice of $\L(P)$ is well-founded and generates $\F(P)$.
\item 
The class $\gwqo^{+}$ contains all the posets $P$ such that $\F(P)$ is 
positively well quasi-ordered generated, 
that is there is a sublattice of $\L(P)$ that generates $\F(P)$ and is a well quasi-ordering.
\item 
The class ${\gbqo^{+}}$ contains all the posets $P$ such that $\F(P)$ is 
positively better generated, 
that is there is a sublattice of $\L(P)$ that generates $\F(P)$ and is a better quasi-order.
\end{enumerate}
We denote by $\WQO$ and $\BQO$ the class of well quasi-orderings and 
better quasi-orderings.
Figure \ref{figure-one} shows known relations between the above classes.
\medskip

More precisely, we have:
\vspace{-2mm}
\begin{itemize}
\item[(1)] 
$P \in \gwf $ iff $P$ is scattered and narrow (part of Theorem~A).
\vspace{-3mm}
\item[(2)]
$\gwqo^+ \subseteq \WQO$ (Theorem~\ref{thm-2.1}).
\vspace{-3mm}
\item[(3)]
$\gwf^+ \supseteq \WQO$ (Theorem~\ref{thm-2.5}).
\vspace{-3mm}
\item[(4)] If $R$ is the Rado poset, then 
$R \in \WQO \setminus \BQO$, and the fact that
$R \in \gwqo^+$ is a direct consequence of Theorem~\ref{thm-2.12}.%
\vspace{-3mm}
\item[(5)] If $P$ is the chain $\omega^* \cdot \omega_1$, then
$P \in \gwf \setminus \gwqo$ (Theorem~\ref{thm-3.1}).%
\vspace{-3mm}
\item[(6)] $\gbqo {\setminus} \WQO \neq \emptyset$ is proved in Theorem~\ref{thm-3.7}.
\end{itemize}
\vspace{-2mm}
\medskip
In section~\ref{section-classw},
we study the class $\calB \eqdef \gbqo^{+}$ .
We show that a member of $\classw$ must be well quasi-ordered,
and that every countable well quasi-ordering is a member of~$\classw$.
We thus ask the following questions.
\begin{question}
\begin{rm}
\label{question-1.2}
(1) Is $\WQO\subseteq \gwqo$\,? 
That is, is it true that every well quasi-ordered set $W$ has a well quasi-ordered generating
lattice for $\F(W)$? 
Similarly, is  $\WQO\subseteq \gbqo$\,?

(2) Is $\WQO\subseteq \gwqo^{+}$\,?  
Is $\WQO\subseteq \gbqo^{+}$\,? 
In other words, is there a well quasi-ordered (better quasi-ordered)
generating lattice for $\F(W)$ contained in $\L(W)$ whenever $W$ is a well quasi-order?

(3) Is there a well founded poset $W$ in $\gwqo\setminus \gbqo^{+}$\,? 
That is, is there a well founded poset $W$  such that
$\F(W)$ has a well quasi-ordered generating lattice,
but $\F(W)$ has no better quasi-ordered generating lattice for $\F(W)$ contained in $\L(W)$?
\end{rm}
\end{question}
Another type of results concerns homomorphic images of poset algebras.
Let $\alpha$ be an ordinal and let $f$ be 
a homomorphism from the ordinal algebra $\F(\alpha)$
 onto a Boolean algebra $A$. 
Then $A$ is easily seen to be an ordinal algebra, and $A$ is generated by $f[\alpha]$ 
and $A$ is isomorphic to $F(f[\alpha])$.
There is no analogous property for poset algebras:
in Theorem~\ref{thm-3.10}, it is shown that there is a well quasi-ordering $W$ 
(which is the ``incomparable sum'' of two copies of~$\omega_1$) 
such that $\F(W)$  has a homomorphic image which is not isomorphic to any poset algebra.

The following theorem~\cite[Theorems~1.4]{ABKR} points 
to some affinity that scattered narrow posets have with well quasi-orders.
\medskip
\newline
{\bf Theorem B}
\begin{it}
Let $P$ be a narrow scattered partially ordered set.
Then there is a well quasi-ordering $W$ and a subalgebra $B$ of $\F(W)$ 
such that $\F(P)$ is a homomorphic image of $B$.
\end{it}
\medskip
\newline
One may ask if Theorem~B can be strengthen by requiring that $B=\F(W)$.
Is it true that every scattered and narrow poset $P$ has a well quasi-ordering $W$ 
such that $\F(P)$ is a homomorphic image of $\F(W)$? 
We provide a negative answer in Theorem~\ref{thm-3.2}:  
we show that there is a scattered and  narrow poset $P$,
in fact the chain $\omega^* \cdot \omega_1$, 
such that $\F(P)$ is not a homomorphic image of $\F(W)$ for any wqo $W$.

This result leaves another possible way of strengthening  Theorem~B: 
to replace the conclusion that $\F(P)$ is a homomorphic image of a subalgebra 
by the stronger statement that $F(P)$ is actually a subalgebra of $\F(W)$. 
We ask the following.
\begin{question}
\label{question-1.3}
\begin{rm}
Let $P$ be a narrow and scattered poset.
Is it true that there is a well quasi-ordering $W$ such that
$\F(P)$ is embeddable in $\F(W)$? 
\end{rm}
\end{question}
The answer to Question~\ref{question-1.3} is positive 
when $P$ is scattered and covered by finitely many chains
(see~\cite{BR2}).

Section~\ref{examples} contains other examples distinguishing 
certain classes of Boolean algebras related to well quasi-orderings.
There is a well generated Boolean algebra which is not embeddable 
in any well quasi-ordered poset algebra (Theorem~\ref{thm-3.4}) and
there is a Boolean algebra which is generated by a better quasi-ordered set 
while it is not isomorphic to any well quasi-ordered poset algebra (Theorem~\ref{thm-3.7}).
\section{Results on the class $\classw$}\label{section-classw}
\label{S2}
We consider in this section the class $\calB \eqdef  \gbqo^{+}$ of 
all posets $P$ for which the algebra $\F(P)$ is generated 
by a better quasi-ordering contained in $\L(P)$.
Clearly, $\calB \subset \gbqo$ and a simple example shows that these two classes are distinct: 
$\omega^{\ast}$ the inverse ordering of $\omega$ is in $\gbqo\setminus \calB$.
Also, all countable scattered chains are in that difference.
Indeed, if $C$ is a countable scattered chain, 
then $\F(C)$ is superatomic and thus isomorphic to $\F(\alpha)$
for some countable ordinal $\alpha$ (see~\cite[\S17.2]{kopp}).

We shall prove that $\classw$ is contained in the class of  well quasi-orderings.
We shall also prove some preservation results on the class~$\classw$
and we show that a minimal well quasi-ordering which 
is not in $\classw$~(assuming it exists) must have uncountable cofinality.

Clearly, by Proposition~\ref{prop-2.3}(c),  all better quasi-orderings belong to $\classw$.
We shall prove next that all posets in $\gwqo^{+}$ are wqo.
\begin{theorem}
\label{thm-2.1}
$\gwqo^{+}\subseteq \WQO$. In words:
If $P$ is a poset such that there exists a well quasi-ordered poset
$Q\subseteq \L(P)$ which generates $\F(P)$.
Then $P$ is a well quasi-ordering.
\end{theorem}
\begin{pf} 
Suppose $P \eqdef \pair{P}{\leq}$ is not well quasi-ordered.
Then there exists a linearly ordered set $R \eqdef \pair{P}{\preceq}$ 
which is a linear augmentation of $P$ 
(that is: for every $p,q \in P$, if $p \leq q$ then $p \preceq q$) and 
which is not well ordered (see~\cite{higman,wolk}).
Since the identity function from $P$ onto $R$ is increasing and onto,
the embeddings $\fnn{x}{P}{\F(P)}$ and $\fnn{x^R}{R}{\F(R)}$
define a homomorphism $\fnn{h}{\F(P)}{\F(R)}$ such that $h(x_p) = x_p^R$.
Since $\setn{ x^R_p }{ p\in P }$ generates~$F(R)$,  $h$ is onto $F(R)$. 
Note that $h[\L(P)] = \L(R)$.
Since $R$ is a chain, $\L(R) = \setn{ x_p^R }{ p \in P }$.
Hence $h[Q] \subseteq \setn{ x_p^R }{ p \in P }$.
If $h[Q] \neq \setn{ x_p^R }{ p \in P }$, 
then the Boolean algebra generated by $h[Q]$ is a proper subalgebra of $\F(R)$.
Hence $h[Q] = \setn{ x_p^R }{ p \in P }$ and thus $h[Q]$ is not well-founded.
On the other hand, since $h$ is increasing and $Q$ is a wqo,  
$h[Q]$ is necessarily also a wqo, which is a contradiction.
\end{pf}
For an exposition of the theory of better quasi-orderings we refer 
to~\cite{fraisse, laver2, milner}.
We will review some definitions and notations.
 A {\it barrier\/} on an infinite subset $S$ of $\omega$
is a family $B$ of pairwise $\subseteq$-incomparable non-empty finite subsets
of $S$ such that for every infinite set $A \subseteq S$
there exists $s \in B$ with $A \cap [0,\max(s)+1) = s$, i.e. $s$ is an initial segment of $A$.
For $s,t \in B$ we write $s \tleq t$ and we say that $s$
{\it precedes\/} $t$ if $s {\setminus} \sngltn{\min(s)}$ is an initial segment of $t$.
In other words, $s \tleq t$ if $s=\setd{n_0}{n_k}$ and $t=\setd{n_1}{n_{k+\ell}}$, 
where $n_0<n_1<\dots<n_{k+\ell}$
(formally it could be that $t=\setd{n_1}{n_k}$ but 
we require that $t\not\subseteq s$ whenever $s,t$ are elements of a barrier).

A poset $P$ is {\it better quasi-ordered} (briefly {\it bqo\/}) if
for any barrier $B$ on $\omega$, for any function $\fnn{f}{B}{P}$ there exist $b_0,b_1\in B$ 
such that $b_0 \tleq b_1$ and $f(b_0)\leq f(b_1)$.
Let $B$ be a barrier such that $\bigcup B = S$.
An infinite subset $B'$ of $B$ is a {\it subbarrier of $B$} if
$B'$ is a barrier on $S' \eqdef \bigcup B'$. 
The Nash-Williams partition theorem says that if $B$ is a barrier 
and $\fnn{f}{ B }{ 2 }$ a partition of $B$ into two classes, 
then there exists a homogeneous sub-barrier of $B$.
The square $B^{2}$ of a barrier $B$ is the collection of all sets 
of the form $s\cup t$ where $s,t\in B$ and $s\tleq t$.
If $\fnn{g}{ B^2 }{ 2 }$ is a partition, 
then some sub-barrier $E$ of $B^2$ is homogeneous and 
it turns out that there exists $B_0$ a sub-barrier of $B$ such that $B_0^2$ is homogeneous.
The following fact about barriers will be used in this work.
Let $P$ be a poset, $B$ a barrier and $\fnn{f}{B}{P}$ be a function.
Then there is a subbarrier $B'$ of $B$ such that:
either 
(1):~$f$ is {\it bad}, that is $f(b_1) \not\leq f(b_2)$ for every $b_1 \tleq b_2$ in $B'$,
or (2):~$f$ is {\it perfect}, that is $f(b_1) \leq f(b_2)$ for every $b_1 \tleq b_2$ in $B'$.

If $X$ is a set and $\kappa$ is a cardinal, 
then $[X]^{<\kappa} \eqdef \setn{ Y \subseteq X }{ |Y|<\kappa }$.

Let $P$ be a poset.
We say that $Q$ is a {\it subposet\/} of $P$ if $Q \subseteq P$ and
$\leq^Q \,\,=\,\, \leq^P\kern-3pt{\restriction}Q$.
For $p \in P$, we set $P^{\geq p} = \setn{ q \in P }{ q \geq p }$,
$P^{\leq p}$, $P^{<p}$ and $P^{>p}$ are defined similarly.
If $A \subseteq P$, 
then $P^{\geq A} \eqdef \bigcup\setn{ P^{\geq p} }{p \in A}$.
The sets $P^{\leq A}$, etc. are defined similarly.

The set $R$ is an {\it initial segment of} $P$
if for every $r \in R$, $P^{\leq r} \subseteq R$.
We denote by $\Is(P)$ the set of all initial segments of $P$, and
by $\Is^{\rm fin}(P)$ the set of all finitely generated initial segments of $P$,
that is, $J \in \Is^{\rm fin}(P)$ iff $J = P^{\leq\sigma}$ for some finite subset $\sigma$ of $P$.
Similarly, we define the notion of a {\it final segment}. The set of all final segments of $P$ will be denoted by $\Fs(P)$.

We summarize some results on well quasi-orderings and 
better quasi-orderings in the following propositions.
\begin{proposition}
\label{prop-2.2}
\rfs{(a)} For every poset the following implications hold: 

\centerline{well ordering $\Rightarrow$
bqo $\Rightarrow$ wqo $\Rightarrow$ well founded.}

\rfs{(b)} Let $P$ be a poset and let $Q$ be a subposet of $P$.
If $P$ is well quasi-ordered (better quasi-ordered), then so is $Q$.

\rfs{(c)}
Let $\pair P\leq$ be a poset and let $\setn{R_i}{i<\ell}$ be a finite
set of ordering relations whose union is $\leq$.
If for every $i<\ell$\,, \,$P_i$ is well quasi-ordered (better quasi-ordered), then so is $P$.

\rfs{(d)} Let $\fnn{f}{P}{Q}$ be an order preserving surjection.
If $P$ is well quasi-ordered (better quasi-ordered), then so is $Q$.

\rfs{(e)}  
If $\pair{P}{\leq}$ is well quasi-ordered (better quasi-ordered),
then so is\\ $\pair{\Is^{\rm fin}(P)}{\subseteq}$.

\rfs{(f)}  
If $\pair{P}{\leq}$ is better quasi-ordered, then $\pair{\Is(P)}{\subseteq}$ is better quasi-ordered.

\rfs{(g)} $\pair{P}{\leq}$ is well quasi-ordered if and only if 
$\pair{\Is(P)}{\subseteq}$ is well-founded.
\end{proposition}
A consequence of Lemma~\ref{prop-2.2} is the following proposition.
\begin{proposition}
\label{prop-2.3}
\rfs{(a)} Let $M$ be a subset of a distributive lattice $L$ such that $M$
is meet-closed and $M$ generates $L$.
\newline
\rfs{(a1)}
Then every member of $L$ is a finite join of members of $M$.
\newline
\rfs{(a2)}
If $M$ is a well founded, then the same holds for~$L$.
\newline
\rfs{(a3)}
If $M$ is a well quasi-ordering (better quasi-ordering), 
then the same holds for~$L$.

\rfs{(b)}
For every poset $P$: if\/ $\Pi(P)$ is well quasi-ordered (better quasi-ordered) 
as a subset of $F(P)$, then so is $\L(P)$.
Recall that  $\Pr(P)$ the meet semilattice of $\F(P)$ generated by
$x[P] \eqdef \setn{ x_p }{ p \in P }$. 
So a member of $\Pr(P)$ is of the form
$x_{\sigma} \eqdef \prod\setn{ x_p }{ p\in\sigma }$
where $\sigma$ is a finite subset of $P$.
 
\rfs{(c)} Let $B$ be a Boolean algebra and let $Q \subseteq B$.
If $Q$ is a better quasi-ordering, 
then so is the lattice generated by $Q$.

\rfs{(d)} If $Q\subseteq \L(P)$ is a better quasi-ordering and generates $\L(P)$ as a lattice, 
then $P$ is a better quasi-ordering.
\end{proposition}

\begin{pf} (a) (a1) is trivial since the lattice is distributive.

(a2) see~\cite[Lemma~2.8(b)]{BR1} and is reported below 
from~\cite[Lemma 3.3(b)]{ABKR}).
For completeness, we recall the proof of~(a2).
It is easy to check that the following holds.
\vspace{-2mm}
\begin{itemize}
\item[$(*)$]
If $w_0>w_1>\cdots $ \
is a strictly decreasing sequence in $L$, and 
$w_0 = \sum_{i<n} v_i$, then there is $\ell <n$ such that 
$\setn{w_j \cdot v_\ell }{ j<\omega }$ contains a strictly decreasing infinite subsequence.
\end{itemize}
\vspace{-2mm}
The proof uses the distributivity of $L$.
Next, suppose by contradiction that $u_0>u_1>\cdots $ \ is a
strictly decreasing sequence in $L$. We define by induction a
strictly decreasing sequence $\seqn{ v_n }{ n \in \omega }$ in $M$. 
Assume by induction that $v_n$ has the following
property. 
There is a strictly decreasing sequence   $w_0>w_1>\cdots $ \ in $L$ such that $w_0<v_n$. 
Let $U\subseteq M$ be a finite set such that $u_0 = \sum U$.
By $(*)$, there is $v_0 \in U$ such that $\setn{u_j\cdot v_0 }{j\in\omega}$ 
contains a strictly decreasing subsequence. 
Hence $v_0$ satisfies the induction hypothesis. 
Suppose that $v_n$ has been defined, and let $v_n > w_0 > w_1 \cdots $ \ 
be as in the induction hypothesis. 
Let $W \subseteq M$ be a finite set such that $\sum W = w_0$. 
By $(*)$, there is $v_{n+1} \in W$ such that $\setn{ w_j \cdot v_{n+1} }{j\in\omega}$ 
contains a strictly decreasing sequence. 
So $v_{n+1}$ satisfies the induction hypothesis and $v_n > w_0 \geq v_{n+1}$. 
The sequence $\seqn{v_n }{n\in\omega} \subseteq M$ is strictly
decreasing. This contradicts the well foundedness of $M$, 
so $T$ is well founded.

(a3).
Suppose that $M$ is wqo (bqo).
By Lemma~\ref{prop-2.2}(e), 
$\pair{[M]^{<\omega}}{\leq}$ is wqo (bqo).
Let $\fnn{f}{[M]^{<\omega}}{L}$ be defined by $f(\tau) = \sum\tau$.
It is easy to check that $f$ is an increasing surjection.
By Lemma~\ref{prop-2.2}(d) $L$ is wqo (bqo).

(b) follows from Part~(a3).

(c)
Let $Q$ be a poset.
We define a partial ordering on the set
$[Q]^{<\omega}$ of finite subsets of~$Q$.
Let $\sigma,\tau \in [Q]^{<\omega}$. Define $\sigma \leq^m \tau$,
if for every $j \in \tau$ there is $i \in \sigma$
such that $j \leq i$.
So $\sigma \leq^m \tau$ iff $Q^{\geq\sigma} \supseteq Q^{\geq\tau}$,
that is $Q-Q^{\geq\sigma} \, \subseteq \, Q-Q^{\geq\tau}$
(as initial segments).
By Lemma~\ref{prop-2.2}(f), $\pair{[Q]^{<\omega}}{\leq^m}$ is a bqo.
Let $\fnn{f}{[Q]^{<\omega}}{B}$ be defined by $f(\sigma) = \prod\sigma$.
Let $M \eqdef \rfs{Rng}(f)$.
Since $f$ is an order preserving function
from $\pair{[Q]^{<\omega}}{\leq^m}$ onto $\pair{M}{\leq}$,
by Lemma~\ref{prop-2.2}(d) $\pair{M}{\leq}$ is bqo.
Next let $\fnn{g}{[M]^{<\omega}}{\L(Q)}$ be defined by
$g(\tau) = \sum\tau$ and let $M' \eqdef \rfs{Rng}(g)$.
Again by Lemma~\ref{prop-2.2}(a), $\pair{M'}{\leq}$ is bqo. 
Now, Part~(c) follows from the fact that $M'$ is the lattice 
generated by $Q$.

(d) Since $Q \subseteq \L(P)$ and $Q$ is a better quasi-ordering,
by Part~(c), the same holds for the sublattice $M$ generated by $Q$.
We have $M \subseteq \L(P)$, and since $Q$ generates $\L(P)$,
$M = \L(P)$ and thus $\L(P)$ is a better quasi-ordered set.
Note that $p \mapsto x_p$ is an embedding from $P$ into $\L(P)$.
By Proposition~\ref{prop-2.2}(b), $P$ is a better  quasi-ordered set.
\end{pf}

The following Theorem~\ref{thm-2.5} is one of the motivations for introducing the class~$\classw$.
Recall that for a finite set $\rho$ of $P$, $x_{\rho}$ denotes $\prod\setn{ x_p }{ p\in\rho }$,
and by definition, $\Pr(P)$ is the set of such $x_{\rho}$'s.
For the proof of Theorem~\ref{thm-2.5}, we use general results on poset algebras
stated in Parts (1) and (2) of Fact~\ref{fact-2.4}.
Parts (3) and (4) of that fact are used to prove Lemma~\ref{lemma-2.11}.
We recall that if $B$ is a Boolean algebra
and $a \in B \setminus \sngltn{0}$,
then $B {\restriction} a$ denotes the Boolean algebra $A$
defined as follows: 
$A \eqdef \setn{ b \in B }{ b \leq a }$,
and the operations are:
$b +^A\, b' = b +^B b'$, $b \cdot^A\, b' = b \cdot^B \, b'$ and
$-^A \, b = a -^B b$. Hence $1^A = a$.
Notice that $A = \setn{ b \cdot a \in B }{ b \in B }$ and
that $\fnn{\varphi}{B}{A}$ defined by $\varphi(b) = b \cdot a$ is a
homomorphism from $B$ onto $A$.

\begin{fact}
\label{fact-2.4}
Let $P$ be a poset.
\begin{enumerate}
\item[\rfs{1.}] Let $\sigma,\tau$ be finite subsets of $P$.
The following properties are equivalent.

\rfs{(i)}
$(\,\prod_{p \in \sigma} x_p)\cdot (\,\prod_{q \in \tau} -x_q) = 0$

\rfs{(ii)}
there are $p \in \sigma$ and $q \in \tau$ such that $p \leq q$.

\item[\rfs{2.}] Let $\sigma,\tau$ be finite subsets of $P$.
(Recall that for a finite subset $\rho$ of $P$, we use the notation $x_{\rho} \eqdef \prod\setn{ x_p }{ p\in\rho }$.)
The following properties are equivalent.

\rfs{(i)}
$x_{\sigma} \leq x_{\tau}$\,.

\rfs{(ii)}
for every $q \in \tau$ there is $p \in \sigma$ such that $p \leq q$.

\rfs{(iii)}
$P^{\geq\sigma} \supseteq P^{\geq\tau}$.

\item[\rfs{3.}]
If $ \fnn{ \imath }{ Q }{ P }$ is an embedding of poset $Q$ into $P$ and 
$\fnn{x^Q}{Q}{\F(Q)}$,   $\fnn{x}{P}{\F(P)}$ are formed,
then there is a Boolean embedding $\fnn{h}{\F(Q)}{\F(P)}$ 
such that $h(x^Q_q) = x_{\imath(q)}$ for all $q\in Q$. 
Moreover  $h[\L(Q)]\subseteq \L(P)$. 

\item[\rfs{4.}]
Let $q \in P$. We set $Q = \setn{ p \in P }{ p \not\geq q }$.
So, $Q$ is a subposet of $P$.
Let $\fnn{ f }{ \F(Q) }{ \F(P) {\restriction} x_q }$ defined by $f(y) = y \cdot x_q$.
Then $f$ is an isomorphism from $\F(Q)$ onto $\F(P) {\restriction} x_q$.
\end{enumerate}
\end{fact}
\begin{proof}
(1) is proved in a direct way in~\cite[Proposition~2.5]{ABKR}.
An alternative proof, obtained via a topological definition of the
Stone space of $\F(P)$ can be found in the introduction of~\S\ref{S3}.

(2) follows from Part~(1).

(3) follows from \cite[Proposition 2.5(d2)]{ABKR}.

(4) The map $\fnn{ h }{ \F(P) }{ \F(P) {\restriction} x_q }$ 
defined by $h(y) = y \cdot x_q$ is a homomorphism onto.
Since $Q \subseteq P$, by Part~(3), we can view $\F(Q)$ as 
a subalgebra of $\F(P)$.
Let $f = h {\restriction} \F(Q)$.
So $f$ is a homomorphism from $\F(Q)$ into $\F(P) {\restriction} x_q$, and
$f(y) = y \cdot x_q$ for every $y \in \F(Q)$.

We show that $f$ is one-to-one.
It suffices to prove that for finite subsets $\sigma$ and $\tau$ of $Q$, setting
$e \eqdef (\,\prod_{s \in \sigma} x_s)\cdot (\,\prod_{t \in \tau} -x_t)$ 
(that is a member of $\F(Q)$):%
\medskip
\newline
\phantom{i}
\hfill
if $f(e) = 0$ then $e = 0$.
\hfill
\phantom{i}
\medskip
\newline
(Recall that $f(e) \eqdef e {\cdot} x_q \in \F(P) {\restriction} x_q$.)
Indeed any member of $\F(Q)$ is a finite sum of such~$e$.
Notice that:
\medskip
\newline
\phantom{i}
\hfill
$f(e) = e \cdot  x_q  
= (\,\prod_{s \in \sigma} x_s) \cdot (\,\prod_{t \in \tau} -x_t)) \cdot x_q
= (\,\prod_{s \in \sigma\cup\sngltn{q}} x_s) \cdot (\,\prod_{t \in \tau} -x_t))$.%
\hfill
\phantom{i}%
\medskip%
\newline%
Suppose $f(e)=0$.
By Part~(2), let $s \in \sigma\cup\sngltn{q}$ and $t \in \tau$
be such that $s \leq t$.
The case $s = q$ does not occur: otherwise $t \geq s = q$,
and thus $t \not\in Q$, that contradicts the fact that $t \in \tau \subseteq Q$.
Next, if $s \in \sigma$, then, by Part~(2) again, $e = 0$.
We have proved that $f$ is one-to-one.%
\smallskip

Next, we show that $f$ is onto.
Since $f$ is a homomorphism, and, by the definition,
since $\F(P)$ is generated by $x[P] = \setn{ x_p }{ p \in P}$,
$\F(P) {\restriction} x_q$ is generated by 
$Y \eqdef \setn{ x_p \cdot x_q }{ p \in P}$.
Note that $x_q \in Y$ and $x_q$ is the unity of $\F(P) {\restriction} x_q$.
So:
\smallskip
\newline
\phantom{i}
\hfill
$\F(P) {\restriction} x_q$ is generated by 
$Y^- \eqdef \setn{ x_p \cdot x_q }{ p \in P \text{ and } x_p \cdot x_q \neq x_q }$.
\hfill
\phantom{i}
\smallskip
\newline
Since $Y^-$ generates $\F(P) {\restriction} x_q$, 
it suffices to prove that for every $y \in Y^-$ there is $x \in \F(Q)$ such that $f(x)=y$.
Let $x_p \cdot x_q \eqdef y \in Y^-$ with $p \in P$.
If $p \not\geq q$ then $p \in Q$ and $f(x_p) = x_p \cdot x_q = y$.
Next, suppose $p \geq q$. Then $y \eqdef x_p \cdot x_q = x_q$ and thus $y \not\in Y^-$:
a contradiction. So $f$ is onto.

We have proved that $f$ is an isomorphism between $\F(Q)$ and
$\F(P) {\restriction} x_q$.%
\end{proof}
\begin{theorem}
\label{thm-2.5}
Let $P$ be a well quasi-ordering.
Then: 

\rfs{1.} $\L(P)$ is well founded,

\rfs{2.} $\L(P)$ is a countable union of well quasi-ordered sets,
and

\rfs{3.} Every antichain of $\L(P)$ is countable.
\end{theorem}
\begin{pf}
(1)
Since $P$ is well quasi-ordered, every final segment $F$ of $P$
is finitely generated, that is, $F$ is of the form
$P^{\geq\sigma}$ where $\sigma \subseteq P$ is finite. 
Hence, by Fact~\ref{fact-2.4}(2), the function $\fnn{f}{ \Is(P) }{ \Pr(P) }$
defined by $f( P\kern1pt{\setminus} P^{\geq\sigma} ) = x_\sigma$
is an order isomorphism between 
$\pair{\Is(P)}{\subseteq}$ and $\pair{\Pr(P)}{\leq}$.
Since $\pair{\Is(P)}{\subseteq}$ is well founded (by  Lemma~\ref{prop-2.2}(g)),
the same holds for $\Pr(P)$.
Since $\L(P)$ is the lattice generated by the meet-closed semilattice $\Pr(P)$,
by Proposition~\ref{prop-2.3}(a2), $\L(P)$ is well founded.

(2) Let $Q_n$ denote the set of all elements of the form 
$x_{\sigma} \in \Pr(P)$ such that $|\sigma|\leq n$.
Since a finite product of wqo sets is wqo, as in the proof of 
Part~(a3) of Proposition~\ref{prop-2.3},
$Q_n$ is wqo.
Now let $L_n(P)$ denote the subset of $\L(P)$ consisting of all
elements of the form $\sum_{i<k} x_{\sigma_i}$, where
$\sigma_i\in [P]^{\leq n}$, $i<k$ and $k<\omega$ is arbitrary.
Then, by Lemma~\ref{prop-2.2}(e), $L_n(P)$ is wqo,
because it is the join sub-semilattice  of $\L(P)$ 
generated by $Q_n$.
Finally, $\L(P) = \bigcup_{n\in\omega} L_n(P)$.

(3) follows from (2).
\end{pf}
\begin{remark}
\label{remark-2.6}
\begin{rm}
In the proof of~\ref{thm-2.5}(1), we have shown 
that $f \,:\, P\kern1pt {\setminus} P^{\geq\sigma} \mapsto x_\sigma$ 
is an order isomorphism between  $\pair{\Is(P)}{\subseteq}$ and $\pair{\Pr(P)}{\leq}$.
If $R$ is the Rado poset, then $\pair{\Is(P)}{\subseteq}$ has an infinite antichain,
and thus the same holds for $\pair{\Pr(P)}{\leq}$.
So $\L(R)$ has an infinite antichain.
\end{rm}
\end{remark}
Below we prove that the class $\classw$ is closed under typical
operations on posets.
\begin{theorem}
\label{thm-2.7}
\rfs{(a)} 
Let $P$ and $Q$ be posets. 
If $\fnn{f}{ P }{ F(Q) }$ is an increasing function such that $f[P]Ê\subseteq L(Q)$
then the extension homomorphism $\fnn{ \hat{f} }{ F(P) }{ F(Q) }$ 
satisfies $\hat{f}[L(P)] \subseteq L(Q)$.

\rfs{(b)} 
Assume $P\in\classw$ and $\fnn{f}{P}{Q}$ is an order preserving
surjection. Then  $Q\in\classw$.
\end{theorem}
\begin{pf} 
(a) follows from the definitions.

(b) 
Let $G \subseteq \L(P)$ be bqo, generating $\F(P)$.
Let $\fnn{\hat{f}}{\F(P)}{\F(Q)}$ be the unique
epimorphism which extends $f$. Then $\hat{f}[\L(P)] = \L(Q)$.
By Lemma~\ref{prop-2.2}(d), $\hat{f}[G]$ is a 
better quasi-ordered subset of $\L(Q)$ which generates $\F(Q)$.
\end{pf}
\begin{theorem}
\label{thm-2.8}
\rfs{(a)} Assume that $Q$ is a better quasi-ordered set and
assume $\setn{P_\xi}{\xi\in Q}\subseteq \classw$.
Then the lexicographic sum
$P= \sum_{\xi\in Q}P_\xi$ is in $\classw$.

\rfs{(b)} Let $(P,\leq^{P})$ be a poset and suppose that
 $\leq^{P}=\leq_0\cup \cdots \cup \leq_{k-1}$\,, 
where each $\leq_{i}$, an ordering of  $P_i$, is in $\classw$.
Then $(P,\leq^{P})\in\classw$.
\end{theorem}
\begin{proof}
(a)
We start with a general result on better quasi-ordered sets.
\begin{itemize}
\item[($*$)] Let $\seqn{ Q_\xi }{\xi \in Q}$ be a family of better quasi-ordered sets, 
indexed by a better quasi-ordered set $Q$.
Then the lexicographic sum $\sum_{\xi \in Q} Q_\xi$
is better quasi-ordered.
\end{itemize}
Fix a barrier $B$ and a function $\fnn{f}{B}{\sum_{ \xi \in Q } Q_\xi}$. 
Consider the square barrier $B^{2}$ and the function $\fnn h{B^2}{2}$ 
defined by $h(s\cup t)= 0$ iff $h(s)$ and $h(t)$ are in the same $P_{\xi}$. 
There is a homogeneous sub-barrier and according to its color we have two possibilities.
(1) For some subbarrier $B'$, $f[B']$ is contained in the same $Q_\zeta$, 
then we can find $b,b'\in B'$ such that $b'$ is a successor of $b$ and $f(b)\leq f(b')$,
because $Q_\zeta$ is bqo.
(2) For  some subbarrier $B'$ of $B$,  $f(b)\in Q_{g(b)}$ and $g(b)<g(b')$ whenever
$b'$ is a successor of $b$ in $B'$.
This shows~($*$).
\smallskip

Since $P_{\xi} \subseteq P$, by Fact~\ref{fact-2.4}(3), we can view $\F(P_{\xi})$ as 
a subalgebra of $\F(P)$.
We assume that $P_\xi \cap P_\zeta=\emptyset$ whenever $\xi \neq \zeta$.
For each $\xi\in Q$, let $G_\xi \subseteq \L(P_\xi)$ be a bqo set generating $F(P_\xi)$.
We may assume that $0 \notin G_\xi$.
Let $G \eqdef \bigcup_{\xi \in Q} G_\xi$.
Clearly, $G \subseteq \L(P)$.

Note that the following holds.
\begin{itemize}
\item[($*$$*$)] For $\xi<\zeta$ in $Q$: 
if $g \in G_\xi$ and $h \in G_\zeta$ then $g < h$.
\end{itemize}
It is clear that $G$ generates $\F(P)$
because $\setn{ x_p }{ p \in \bigcup_{ \xi \in Q } P_\xi }$ does.
By ($*$), $\sum_{\xi \in Q} G_\xi$ is a bqo.
Now, $G$ is bqo because by ($*$$*$),
it is an augmentation of the lexicographic sum $\sum_{\xi \in Q} G_\xi$ 
(that is if $h \leq g$ in $\sum_{\xi \in Q} G_\xi$ then $h \leq g$ in $G$)
and $Q$ is bqo.

(b) Let $R$ be the disjoint sum of $P_0, \dots, P_{k-1}$, 
the order of $R$ is $p<p'$ iff $p,p'\in P_i$ 
for some $i$ and $p <^{P_i}  p'$ in $P_i$.
Then $R$ is the lexicographic sum of $\seqn{P_i}{i \in k}$
where $k$ is endowed with the discrete order.
By Part~(b), $R \in \classw$. 
Now observe that the map $\fnn{f}{R}{P}$,
such that $f {\restriction} P_i = \rfs{Id}_{P_i}$, 
is order preserving onto $P$.
By Theorem~\ref{thm-2.7}(b), $P\in \classw$.
\end{proof}
To prove that if $P,Q \in \classw$ then $P{\times}Q \in \classw$
(Theorem~\ref{thm-2.10}), we need a preliminary result.
Let $P$ and $Q$ be posets. 
We form $\F(P)$ and $\F(Q)$ with the associated maps $\fnn{x}{P}{ \F(P)}$ and 
$\fnn{x}{Q}{\F(Q)}$ 
(we use the same letter to denote both functions, which should not creat any confusion).
Recall that for $u \in P$, $x_u \in x[P]$. 
Form also $\F(P\times Q)$ with the associated embedding 
$\fnn{x}{ P{\times}Q }{  \F(P{\times} Q) }$.

\begin{lemma}
\label{lemma-2.9}
Let $P$ and $Q$ be posets. 
There exists a function $\fnn{E}{ \L(P){\times}\L(Q) }{ \L(P{\times}Q) }$ 
such that the following hold:
\begin{enumerate}
\item[\rfs{1.}]
For every $p\in P$ and $q\in Q$, $E(x_{p},x_{q}) = x_{(p,q)}$.
\item[\rfs{2.}]
For every $a\in \L(P)$ the function taking $b\in \L(Q)$ to $E(a,b)$ can be extended
to a homomorphism from $\F(Q)$ to $\F(P\times Q)$.
\item[\rfs{3.}] For every $x_{q} \in x[Q]$, the function taking $a \in \L(P)$
to $E(a,x_{q})$ can be extended to a homomorphism from $\F(P)$ to $\F(P\times Q)$.
\item[\rfs{4.}]
For every $b\in \L(Q)$, the function taking $a\in \L(P)$ to $E(a,b)$ is order-preserving.
\end{enumerate}
\end{lemma}
\begin{proof} 
For every $q\in Q$ define $\fnn{f_q}{ P }{ \L(P{\times}Q) }$ by the equation 
$f_q(p)=x_{\pair{p}{q}}$. 
Clearly, $f_q$ is order-preserving, and
\begin{equation}
\label{Eq:1}
\text{if}\ q_1< q_2\ (\text{in} \ Q),\ \text{ then: } f_{q_1}(p)< f_{q_2}(p)
\text{\ for any\ } p\in P.
\end{equation}
There is a homomorphism $\fnn{ \hat{f}_q }{ \F(P) }{ \F(P\times Q) }$ that extends $f_q$ 
in the sense that for every $p\in P$,
\[ \hat{f}_q(x_{p})=f_q(p).\]
Since $f_q(p)=x_{(p,q)} \in \L(P{\times}Q$, by Theorem~\ref{thm-2.7}(a),
\begin{equation}
\label{Eq:2}
\hat{f}_q[L(P)] \subseteq \L(P\times Q).
\end{equation}
It follows immediately that
\begin{equation}
\label{Eq:3}
\text{if}\ q_1\leq q_2\ \text{in } Q,\ \text{then for every}\ a\in \L(P)\ 
\hat{f}_{q_1}(a)\leq  \hat{f}_{q_2}(a).
\end{equation}
Fixing $a \in \L(P)$ define $\fnn{g_a}{ Q }{ \F(P{\times}Q) }$ by 
the equation $g_a(q)=\hat{f}_q(a)$.

It follows from (\ref{Eq:3}) that $g_a$ is order-preserving, and by 
(\ref{Eq:2}), 
$g_a[Q]\subseteq \L(P\times Q)$.

The universality property for $Q$ yields a homomorphism 
$\fnn{ \hat{g}_a }{ \F(Q) }{ \F(P{\times}Q) }$ 
that extends $g_a$ in the sense that for $q\in Q$, 
$$\hat{g}_a(x_{q})=g_a(q)=\hat{f}_q(a).$$

We now define for every $a\in \L(P)$ and $b \in \L(Q)$, 
$$E(a,b)=\hat{g}_a(b).$$

It is clear from the definition that for a fixed $a \in \L(P)$ the function 
taking $b\in \L(Q)$ to $E(a,b)$ is extended by the homomorphism $\hat{g}_a$.

Now, for $x_{q}\in \L(Q)$ and $a\in \L(P)$, $E(a,x_{q})=\hat{g}_a(q)=\hat{f}_q(a)$.
So that $\hat{f}_q$ is a homomorphism defined on $\F(P)$ that extends $E(a,x_{q})$ 
(as a function of $a$).

Consider the set $B$ of all $b\in \L(Q)$ such that the function taking $a\in \L(P)$ to $E(a,b)$
is order-preserving.
We have:
\begin{enumerate}
\item $x[Q] \subseteq B$, and
\item If $b_1,b_2\in B$ then $b_1\cdot b_2\in B$ and $b_1+ b_2\in B$
\end{enumerate}
These two items imply that $B=L(Q)$ as required.
The first item follows since for  $x_{q}\in x[Q]$,
if $a_1< a_2$ in $L(P)$ then we have
\smallskip
\newline
\phantom{i}
\hfill
$E(a_1,x_{q}) = \hat{g}_{a_1}(x_{q})= x_{{f}_q(a_1)} \leq x_{{f}_q(a_2)}=E(a_2,x_{q})$.
\hfill
\phantom{i}
\smallskip
\newline
For item (2) notice first that
\smallskip
\newline
\phantom{i}
\hfill
$E(a,b_1\cdot b_2) \eqdef \hat{g}_a(b_1 \cdot b_2) = \hat{g}_a(b_1) \cdot \hat{g}_a(b_2)
= E(a,b_1) \cdot E(a, b_2)$. 
\hfill
\phantom{i}
\smallskip
\newline
Hence if $E(a,b_1)$ and $E(a, b_2)$ are increasing  as functions of $a$, 
the same holds for $E(a,b_1\cdot b_2)$.
So  $b_1 \cdot b_2  \in B$.
Similarly we get that $b_1+ b_2  \in B$. 
This proves the theorem.
\end{proof}
\begin{theorem} 
\label{thm-2.10}
If $P$ and $Q$ are in $\calB$ then so is their product $P\times Q$. 
\end{theorem}
\begin{proof}
Let $A\subseteq \L(P)$ and $B\subseteq \L(Q)$ be better quasi-orderd sets that 
generate $\F(P)$ and $\F(Q)$ (respectively).
Let $C=E[A{\times}B]$ be the image of their product under~$E$.
We claim that $C$ is a bqo subset of 
$L(P{\times}Q)$ that generates $\F(P{\times}Q)$.

Firstly, since $\fnn { E }{ \L(P){\times}\L(Q) }{ \L(P{\times}Q) }$ is order-preserving, 
and since $A{\times}B$ is bqo as product of two bqo's, 
it follows that $E[A{\times}B]$ is bqo.

Secondly, it suffices to prove that for every $p\in P$ and $q\in Q$, 
$x_{\pair{p}{q}}$ is generated by $C$ 
(because $x[P{\times}Q]$ generates $\F(P{\times}Q)$\,). 
Let $C^\ast$ be the Boolean algebra generated by $C$ in $\F(P{\times}Q)$.
\medskip
\newline
{\bf Claim 1.}
\begin{it}
For every $a\in A$, $E(a,x_{p})\in C^\ast$.
\end{it}
\medskip
\newline
\Proof
In $\F(Q)$ we find that $x_{q}= \sum_{i \in I} \prod \beta_i$
where $I$ is finite and $\beta_i$ is a finite sequence of members of $B$ and their complements.
Applying $\hat{g}_a$ which is a homomorphism defined over $\F(Q)$,
\smallskip
\newline
\phantom{i}
\hfill
$\hat{g}_a(x_{q})=\sum_{i\in I} \prod \hat{g}_a[\beta_i].$
\hfill
\phantom{i}
\smallskip
\newline
Here $\hat{g}_a[\beta_i]$ is the sequence obtained from $\beta_i$ by applying $\hat{g}_a$
to each of its members. 
Every member of $\hat{g}_a[\beta_i]$ is of the form $\hat{g}_a(b)$ or 
$-\hat{g}_a(b)$ for some $b\in B$. 
This shows that $E(a,x_{q})=\hat{g}_a(x_{q})$ is in $C^\ast$.
We have proved Claim~1.%
\medskip
\newline
{\bf Claim 2.}
\begin{it}
$E(x_{p},x_{q})=x(p,q)$ is in $C^\ast$.
\end{it}
\medskip
\newline
\Proof
In $\F(P)$ we can present $x_{p} = \sum_{i \in I} \prod \alpha_i$ where 
$I$ is finite and $\alpha_i$ this time is a finite sequence of members of $A$ 
and their complements.
Since $\hat{g}_a(x_{q}) = \hat{f}_q(a)$, and as every member of the form $\hat{f}_q(a)$ 
is in $C^\ast$, we apply 
the homomorphism $\hat{f}_q$ and obtain as before that $\hat{f}_q(x_{p})$
is in $C^\ast$. 
Namely, $x_{\pair{p}{q}} \in C^\ast$ as required.
So Claim~2 is proved.
\smallskip

Since $C \eqdef E[A{\times}B]$ is a bqo subset of  $L(P{\times}Q)$,  
and for every $p\in P$ and $q\in Q$, $x_{\pair{p}{q}}$ is generated by $C$,
$C$ generates $\F(P{\times}Q)$, that is $P{\times}Q \in \classw$.
\end{proof}
Finally, we are going to show that all countable well quasi-orderings
are in $\classw$.
\begin{lemma}
\label{lemma-2.11}
Let $P$ be a directed well quasi-ordered poset with
countable cofinality and such that every proper initial segment
of $P$ is in $\classw$.  Then $P\in\classw$.
\end{lemma}
\begin{pf}
First suppose that $P$ has a last element, that we denote by $p$.
Let $Q = P \setminus \sngltn{p}$.
So $Q$ is a proper initial segment of $P$.
We consider $Q$ as a subposet of $P$.
By the hypothesis, $Q \in \classw$.
Trivially, $\sngltn{p} \in \classw$.
Let $P'$ be the disjoint union of $Q$ and $\sngltn{p}$.
That is $p <^{P'} q$ iff $p,q \in Q$ and $p <^Q q$.
By Theorem~\ref{thm-2.7}(a), $P' \in \classw$.
Since $P$ is the lexicographic sum $Q + \sngltn{p}$
and the identity function from $P'$ onto $P$ is increasing,
by Theorem~\ref{thm-2.7}(b), $P \in \classw$.

Next suppose that $P$ has no last element.
Let $\setn{p(n)}{n\in\natN}$ be a strictly increasing cofinal
sequence in $P$ and define $P_n \eqdef \setn{x\in P}{ x \not\geq p(n) }$.

Let $n\in\natN$.
Since $P_n$ is a subposet of $P$, by Fact~\ref{fact-2.4}(3)
$\F(P_n)$ is a subalgebra of $\F(P)$.
Since $P_n$ is a proper initial segment of $P$\,,
\,$\F(P_n)$ is generated by a better quasi-ordered subset $G_n$ contained in $\L(P_n)$.
Let $G'_{n+1} = G_{n+1} \cup \sngltn{x_{p(n)}}$.
If $x_{p(n)} \in G_{n+1}$ then $G'_{n+1} = G_{n+1}$ and thus $G'_{n+1}$ is bqo.
Suppose that $x_{p(n)} \not\in G_{n+1}$.
Note that $p(n) \in P_{n+1}$, and by Proposition~\ref{prop-2.2}(c),
$G'_{n+1} \eqdef G_{n+1} \cup \sngltn{x_{p(n)}}$ is a bqo generating set for $\F(P_{n+1})$
and $G'_{n+1} \subseteq \L(P_{n+1})$.
So we may assume that: 
\begin{itemize}
\item[{}]
$\F(P_n)$ is generated by a better quasi-ordered subset $G_n \subseteq \L(P_n)$
\\
and $x_{p(n)} \in G_{n+1}$.
\end{itemize}
Since $P_n$ is a subposet of $P$, $\L(P_n) \subseteq \L(P)$.
Let
\medskip
\newline
\phantom{i}
\hfill
$H_n = \setn{ x_{p(n-1)} + y \cdot x_{p(n)} }{ y \in G_n }$
where $x_{p(-1)} \eqdef 0$.
\hfill
\phantom{i}
\medskip
\newline
Let $H \eqdef \sngltn{0} \cup \bigcup_{n\in\omega} H_n$.
Note that $\setn{ x_{p(n)} }{ n\in\omega \cup \sngltn{-1} } \subseteq H $.
\smallskip

We show that $H$ is as required.
That is $H \subseteq \L(P)$, 
$H$ is better quasi-ordering and $H$ generates $\F(P)$. 
\smallskip

Since $G_n \subseteq \L(P_n) \subseteq \L(P)$, 
$H_n \subseteq \L(P)$ and thus $H \subseteq \L(P)$.
\smallskip

Next we show that $H$ is a better quasi-ordering.
Let $\fnn{f}{G_n}{\F(P)}$ defined by $f(y) = x_{p(n-1)} + y \cdot x_{p(n)}$.
By the definition of $H_n$, $\rfs{Rng}(f) = H_n$.
Note that $y \mapsto x_{p(n-1)} + y \cdot x_{p(n)}$
is an order preserving map from the better quasi-ordering $G_n$ onto $H_n$. 
By Proposition~\ref{prop-2.2}(d), $H_n$ is better quasi-ordered.
Next, it is easy to check that:
for $h_m\in H_m$ and $h_n \in H_n$: if $m<n$ then 
$h_m \leq x_{p(m)} \leq x_{p(n-1)} \leq h_n$.
So $H' \eqdef \bigcup_{n\in\omega} H_n$ is order-isomorphic to the lexicographic sum 
$\sum_{n\in\omega} H_n$, and since each $H_n$ is a bqo,  
$H'$ is better quasi-ordered.
(The fact that $\sum_{n\in\omega} H_n$ is bqo follows from the item~($*$) of the proof
of Theorem~\ref{thm-2.8}(a).)
Now, $H$ is the union of the two bqo $H'$ and~$\sngltn{0}$,
and thus $H$ is a better quasi-ordered set.
\smallskip

It remains to show that $H$ generates $\F(P)$.
We show first that if $p \in P$, then $x_p$ is in the Boolean algebra
generated by $H$.

For a subset $X$ of $\F(P)$, we denote by $\rfs{cl}(X)$ the subalgebra of $\F(P)$
generated by $X$.
We begin by a simple observation.
\vspace{-1.5ex}
\begin{eqnarray}
\text{If}\  y \in G_n \  \text{then} \ 
y \cdot (x_{p(n)} - x_{p(n-1)}) 
&
= 
&
(x_{p(n-1)} + y \cdot x_{p(n)}) - x_{p(n-1)}
\nonumber
\\
&
=
&
 f(y)  - x_{p(n-1)}.
\nonumber
\vspace{-1.5ex}
\end{eqnarray}
Since $y \in G_n$, $f(y) \in H_n \subseteq H$, and since $x_{p(n-1)} \in H$,
$f(y)  - x_{p(n-1)} \in \rfs{cl}(H)$. We have proved:
\begin{itemize}
\item[(1)]  if  $y \in G_n$ then 
$y \cdot (x_{p(n)} - x_{p(n-1)}) = f(y)  - x_{p(n-1)} \in \rfs{cl}(H)$.
\end{itemize}
By Fact~\ref{fact-2.4}(4), 
$\fnn{\psi}{ \F(P_n) }{ F(P) {\restriction} x_{p(n)} }$ defined by $\psi(y) = y \cdot x_{p(n)}$
is a homomorphism onto, and obviously 
$\fnn{ \varphi }{ \F(P) {\restriction} x_{p(n)} }{ \F(P) {\restriction} (x_{p(n)} {-} x_{p(n-1)}) }$
defined by $\varphi(z) \eqdef z {\cdot} (x_{p(n)} {-} x_{p(n-1)}) = z  - x_{p(n-1)}$ 
is also a homomorphism onto.
Let $\fnn{g}{ \F(P_n) }{ \F(P) {\restriction} (x_{p(n)} {-} x_{p(n-1)}) }$
defined by $g(z) = z \cdot (x_{p(n)} {-} x_{p(n-1)})$.
So:
\begin{itemize}
\item[(2)]  $g = \varphi \circ \psi$ is a homomorphism onto.
\end{itemize}
Recall that $H_n \eqdef \setn{ f(y) }{ y \in G_n }$. Hence by (1), it follows that
\begin{itemize}
\item[(3)] $g[G_n] = \setn{ t - x_{p(n-1)} }{ t \in H_n }$.%
\end{itemize}
Since $x_{p(n-1)} \in H$ and $H_n \subseteq H$, 
by (3), 
\begin{itemize}
\item[(4)] $g[G_n] \subseteq \rfs{cl}(H)$.
\end{itemize}
By the definition of $g$, and by (2), since $g$ is onto, and since,
by the definition,  $\rfs{cl}(G_n) = \F(P_n)$, 
\begin{itemize}
\item[(5)] $g[G_n]$ generates 
$\F(P) {\restriction} (x_{p(n)} {-} x_{p(n-1)})$.
\end{itemize}
>From (4) and (5), it follows that:
\begin{itemize}
\item[(6)]
$\F(P) {\restriction} (x_{p(n)} {-} x_{p(n-1)}) \subseteq \rfs{cl}(H)$.
\end{itemize}
Next, let $p \in P$.
Recall that $\setn{p(n)}{n\in\natN}$ is strictly increasing and cofinal in $P$.
Let $\ell\in\omega$ be such that $x_p \leq x_{p(\ell)}$.
For $i \leq \ell$, let $y_i \eqdef x_p \cdot (x_{p(i)} {-} x_{p(i-1)})$.
So $y_i \in F(P) {\restriction} (x_{p(i)} {-} x_{p(i-1)})$, and thus,
by (6), $y_i \in \rfs{cl}(H)$.
>From the facts
that $x_p = \sum_{i\leq\ell} \, y_i$, and that $y_i \in \rfs{cl}(H)$ for every $i \leq \ell$,
it follows that $x_p  \in \rfs{cl}(H)$. 

Since $x[P] \eqdef \setn{ x_p }{p \in P}$ generates $\F(P)$,
and $x_p  \in \rfs{cl}(H)$ for any $p \in P$,
$\F(P) = \rfs{cl}(H)$.
\end{pf}
\begin{theorem}
\label{thm-2.12}
Every countable well quasi-ordered poset is in $\classw$.
\end{theorem}
\begin{pf} Suppose not and 
let $P$ be a  countable poset not in $\classw$. Since the set of initial segments 
of a wqo set is well-founded, 
we may assume that every proper segment of $P$ is in $\classw$.
It is well-known (see for instance~\cite[\S4.7.1]{fraisse}), that
every wqo set is a finite union of directed initial
segments (where $X$ is directed if every two members of $X$ have an upper bound in $X$). 
Say $P=P_{0}\cup \cdots \cup  P_{k-1}$. If all $P_{i}$ are proper subsets of $P$, 
then they are all in $\classw$, and it follows from Theorem~\ref{thm-2.8}(b) 
that $P$ itself is in $\calB$.
Otherwise, $P$ itself is directed, and then Lemma \ref{lemma-2.11} applies and 
yields that $P\in \classw$. 
\end{pf}
\section{Examples}
\label{examples}
\label{S3}
In this section we present some examples distinguishing classes
of Boolean algebras related to well quasi-orderings.
The first example, as announced in the introduction, shows that
Theorem B cannot be strengthened by saying 
``every free algebra over a scattered and narrow poset is a
homomorphic image of a wqo poset Boolean algebra''.
The second example is a well generated Boolean algebra which is not
embeddable into a poset algebra and is not a homomorphic image of 
any scattered and narrow poset algebra.
The third one is a scattered and narrow poset algebra which is
generated by a better quasi-ordering while 
not isomorphic to any well quasi-ordered poset algebra.
The last one is a bqo generated Boolean algebra which is not
isomorphic to any poset Boolean algebra.

We need another alternative description of $\F(P)$ using topology.
Viewing a final segment as its characteristic function,
it turns out that $\Fs(P)$, the set of all final segments of $P$,
is a closed subspace of $\{0,1\}^{P}$.
To see this, let $R \subs {P}$, $R\notin\Fs(P)$.
Choose $r,s \in P$ with $r\leq s$, $r\in R$ and $s \notin R$.
Then 
$\setn{ S \subs P }{ r \in S \mbox{\rm\ and } s \not\in S }$
is an open set in $\{0,1\}^{P}$, 
containing $R$, and disjoint from $\Fs(P)$.
We denote by $\widehat{F}(P)$ the Boolean algebra of 
clopen subsets of $\Fs(P)$.
Consider the map
$x_p \mapsto V_p \eqdef \setn{ R \in \Fs(P) }{ p \in R }$.
In \cite[Theorem~2.2]{ABKR} it has been shown that this map extends to 
an isomorphism between $\F(P)$ and $\widehat{F}(P)$.
Thus we have a description of $\F(P)$ as the clopen algebra of the compact space $\Fs(P)$.
It follows that the following holds.
\vspace{-1mm}
\begin{itemize}
\item[($\dagger$)]
Let $\sigma,\tau$ be finite subsets of $P$.
The following properties are equivalent.%
\medskip
\newline 
\rfs{(i)}
$(\,\prod_{p \in \sigma} x_p)\cdot (\,\prod_{q \in \tau} -x_q) = 0$.
\medskip
\newline 
\rfs{(ii)}
$(\,\bigcap_{p \in \sigma} V_p) \cap (\, \bigcap_{q \in \tau} -V_q) = \emptyset$.
\medskip
\newline 
\rfs{(iii)}
there are $p \in \sigma$ and $q \in \tau$ such that $p \leq q$.
\end{itemize}
\vspace{-1mm}
So we have seen that the space $\Fs(P)$ is (homeomorphic to) 
the Stone space $\rfs{Ult}(\F(P))$ of ultrafilters of~$\Fs(P)$.
\medskip

We would like to know which Boolean spaces are homeomorphic to a space of the form $\rfs{Ult}(\F(P))$ for some poset $P$. It turns out that these are precisely those spaces which have a structure of a topological distributive lattice.
Recall that a {\it compact $0$-dimensional lattice $\sixtpl{L}{\tau}{0^L}{1^L}{\meet}{\join}$}
is a compact $0$-di\-men\-sio\-nal space $\pair{L}{\tau}$ such that
$\fifthtpl{L}{0^L}{1^L}{\meet}{\join}$ is a distributive lattice 
with $0^L=\min(L)$, $1^L=\max(L)$ and
the operations $\meet$ and $\join$ are continuous.

If $P$ is a poset, then $\pair{\Fs(P)}{\subseteq}$ is a  compact $0$-di\-men\-sio\-nal lattice.
Conversely, if 
$\sixtpl{L}{\tau}{0^L}{1^L}{\meet}{\join}$ is a 
compact $0$-dimensional lattice, then there is a poset $P$
such that $\pair{L}{\leq} \cong \pair{\Fs(P)}{\subseteq}$. See~\cite[Theorem~2.6]{ABKR} for the details.
\subsection{A superatomic interval algebra which is not a
homomorphic image of a well quasi-ordered poset algebra}
\label{S3a}
This part concerns the fact that  $\gwf {\setminus \gwqo} \neq \emptyset$
and the fact that we cannot strengthen Theorem~B 
by requiring that $B = \F(W)$. 
The example is $C=\omega^*\cdot \omega_1$ which will be
represented as $\omega{\times}\omega_1$ with the order defined by
$\pair{n}{\alpha} <\pair{m}{\beta}$ iff
either $\alpha<\beta$ or $\alpha=\beta$ and $n>m$.
Since $C$ is a scattered chain, $\F(C)$ is well generated,
by Theorem A.
However, it is easy to describe a well founded lattice which
generates $\F(C)$: denote by $G_0$ the set of all atoms of $\F(C)$
and by $G_1$ the set of all $x_{\pair{1}{\alpha}}$
for $\alpha<\omega_1$ (so $G_1$ is of order-type $\omega_1$).
Then the lattice generated by $G_0\cup G_1$
is well founded and generates $\F(C)$.
\begin{theorem}\label{thm-3.1}
Assume that $G$ is a well founded lattice which generates
\break
$\F(\omega^*\cdot \omega_1)$.
Then $G$ contains an uncountable antichain.
\end{theorem}
\begin{pf}
Let $C \eqdef \omega^*\cdot \omega_1$.
We say that $a \subseteq C$ is unbounded in $C$ 
if there is $c \in C$ such that $C^{\geq c} \subseteq a$.
Denote by $\varrho(a)$ the minimal $\alpha<\omega_1$ such that
$C^{\geq\pair{0}{\alpha}} \subseteq a$.
Since $\F(C)$ is an interval algebra, it can be uniquely identified with a subalgebra of $\calP(C)$. Thus we treat elements of $\F(C)$ as subsets of $C$.

Notice that there are only 
countably many unbounded elements $a\in \F(C)$ with a fixed $\varrho(a)$.
Now suppose that  
$\setn{a_n}{n\in\omega} \subseteq G$ is a sequence 
such that each $a_n$ is unbounded and
$\setn{\varrho(a_n)}{n\in\natN}$
is strictly increasing.
Then 
$a_0 \supset a_0 {\cap} a_1 \supset a_0 {\cap} a_1 {\cap} a_2 
\supset \cdots$ 
is a stricly decreasing sequence of members of $G$, 
which contradicts the fact that $G$ is a well founded lattice.
It follows that:
\vspace{-2mm}
\begin{itemize}
\item[(i)] the set $U^G$ of unbounded elements of $G$ in $C$ is countable.
\end{itemize}
\vspace{-2mm}
Let $\alpha \in \omega_1$.
We say that an element $a\in G$ is $\alpha$-good if
$(\omega{\times}\sngltn{\alpha}) \cap a \neq \emptyset$ and
$((\omega \setminus [0,k)) {\times} \sngltn{\alpha})\cap a = \emptyset$ 
for some $k<\omega$.
We claim that
\vspace{-2mm}
\begin{itemize}
\item[(ii)]
there exists an $\alpha$-good element in $G$.
\end{itemize}
\vspace{-2mm}
Suppose otherwise and fix $\alpha\in\omega_1$. Let
\smallskip
\newline
\phantom{abk}\hfill
$H \eqdef
\setn{ \, a\in G }
{ (\exists k \in \omega)
( \,(\omega \setminus [0,k)) {\times} \sngltn{\alpha}\subseteq a \,) \,}$.
\hfill\phantom{abk}
\smallskip
\newline
Obviously $H$ is a sublattice of $G$.
Next, since~$G$ is a set of generators for~$\F(C)$,
($\dagger$):~$G$ separates distinct atoms of~$\F(C)$.
That is: for distinct atoms $a_0$ and $a_1$ of $\F(C)$, 
there are $g \in G$ and $i<2$ such that 
$a_i \subseteq g$ and $a_{1-i} \cap g = \emptyset$.

In particular $G$ separates distinct atoms of $\F(C)$ 
contained in $\omega{\times}\sngltn{\alpha}$,
and therefore $H \neq \emptyset$.
Let $h_0=\min(H)$ 
(which exists because $H$ is a nonempty well founded lattice).
Since $h_0$ is not $\alpha$-good and
$(\omega{\times}\sngltn{\alpha}) \cap h_0 \neq \emptyset$,
we have that 
$(\omega \setminus [0,k)) {\times} \sngltn{\alpha}) \subseteq h_0$
for some $k<\omega$. 
But now, we cannot separate the atoms of $\F(C)$ contained in
$(\omega \setminus [0,k)) {\times} \sngltn{\alpha}$ 
by using elements of $G$,
which is a contradiction.
We have proved~(ii).
\smallskip

For $\alpha\geq1$, choose an $\alpha$-good element $a_\alpha\in G$ and let $f(\alpha)<\alpha$ be such that 
$a_\alpha\cap
([0,\alpha){\times}\omega) \subseteq [0,f(\alpha){+}1){\times}\omega$.
Then the function 
$\fnn{f}{\omega_1 \setminus \sngltn{0}}{\omega_1}$
is regressive, so by Fodor's Theorem, there are $\xi<\omega_1$
and a stationary set $S \subseteq \omega_1$ such that
($\ddagger$):~$a_\alpha\cap([0,\alpha) {\times}\omega) \subseteq [0,\xi{+}1){\times}\omega$
whenever $\alpha\in S$.

Since, by~(i), $U^G$ is countable, let $\alpha_0<\omega_1$ be 
such that $\alpha_0 > \varrho(a)$ for every $a\in U^G$.
By~(i) again, 
\smallskip
\newline
\phantom{i}
\hfill
$T' \eqdef \setn{ \alpha \in S }
{ \alpha > \max(\dbltn{\alpha_0}{\xi+1}) \mbox{\rm\ and } a_\alpha \in G {\setminus} U^G }$ 
\hfill
\phantom{i}
\smallskip
\newline
is uncountable. 
Since for every $\alpha \in T'$ there is $\beta$ such that 
$a_{\alpha} \subseteq \omega{\times}\beta$,
let~$T$ be an uncountable subset of $T'$ 
such that for every $\alpha,\beta \in T$:
(1)~for every $\alpha \in T$, $\alpha>\alpha_0$ and $\alpha>\xi$, and
(2)~if $\alpha<\beta$ then 
$a_{\alpha} \subseteq \omega{\times}\beta$.

Now if $\alpha,\beta\in T$  and $\alpha < \beta$, then,
by the goodness of $a_\alpha$, 
$(\omega{\times}\sngltn{\alpha}) \cap a_\alpha\neq\emptyset$, and,
by~($\ddagger$) and the fact that $\alpha > \xi+1$,
$(\omega{\times}\sngltn{\alpha}) \cap a_\beta=\emptyset$. Thus
$\setn{a_\alpha}{\alpha\in T}$ is an uncountable  antichain in $G$.
\end{pf}
Theorem~\ref{thm-3.1} shows that
$\omega^*\cdot \omega_1 \in \gwf \setminus \gwqo$.
\begin{theorem} \label{thm-3.2}
$\F(\omega^*\cdot\omega_1)$ is not a homomorphic image
of $\F(W)$ for any well quasi-ordered poset $W$.
\end{theorem}

\begin{pf} Suppose $W$ is wqo and
$\fnn{f}{\F(W)}{\F(\omega^*\cdot\omega_1)}$ is an epimorphism.
Since $f$ preserves the lattice operations,
by Proposition~\ref{prop-2.3}(b), $f[\L(W)]$ is a well founded lattice
which generates $\F(\omega^*\cdot\omega_1)$.
By Theorem~\ref{thm-2.5}(3), 
every antichain of $\L(W)$ is countable. 
Since $f$ is increasing, the same holds for $f[\L(W)]$.
This contradicts Theorem~\ref{thm-3.1}.
\end{pf}
Theorem~\ref{thm-3.2} means that the interval space $\Fs(\omega^*\cdot\omega_1)$
(\,$\cong 1+\omega^*\cdot\omega_1 +1$\,)
is not topologically embeddable in the space $\Fs(W)$ for any
well quasi-ordered set~$W$.
\bigskip

The next section is motivated by two facts. 
First every free Boolean algbera is  a poset algebra (over an antichain), 
and every Boolean algebra is a homomorphic image of a free Boolean algebra.
Next every well-generated Boolean algebra $B$ is a homomorphic image of
a free Boolean algebra over a well quasi-ordered poset.
To see this, let $G$ be a well-founded sublattice of $B$ generating $B$.
Then the inclusion mapping $G \subseteq B$ is extendable in a
homomorphism $f$ from $\F(G)$ onto $B$.
We will see that we cannot have this kind of results in the class of poset algebras.
\subsection{A well generated Boolean algebra which is not embeddable into any
poset algebra and which is not an image of any narrow poset algebra}
\label{S3b}
Let $\calA \subseteq \calP(\omega)$ be an uncountable almost disjoint
family and let $B(\calA)$ denote the subalgebra of $\calP(\omega)$
generated by 
$\calG =\setn{ \sngltn{n} }{ n\in\omega } \cup\calA$.
This is called an {\it almost disjoint algebra} and
its topological version is called {\it Mr\'owka's $\Psi$-space}.

In order to show that $B(\calA)$ has the desired properties,
we need a result from the theory of supercompact topological spaces.
A topological space $X$ is {\it supercompact} \cite{vanMill}
if it has a subbase $\calS$ for the closed sets which is 
{\it binary}, i.e. whenever $\calS' \subseteq \calS$ has an 
empty intersection, then  $A_0 \cap A_1 =\emptyset$ for some
$A_0,A_1 \in \calS'$.

If $P$ is a poset, then, by ($\dagger$),  the family 
$\calS = \setn{ V_p }{ p\in P } \cup \setn{ -V_p }{ p\in P }$ is a binary subbase of $\Fs(P)$.
In other words, the Stone space of every poset Boolean algebra is supercompact.

We denote by $A(\kappa)$ the one-point compactification of the discrete space of cardinality $\kappa$.
So $A(\kappa) = \kappa \cup \sngltn{\infty}$.
Note that the clopen algebra of this space is the algebra 
of finite and cofinite subsets of $\kappa$.

A result of Bell \cite[Corollary 3.2]{bell} says that if a compactification $K$ 
of the natural numbers $\natN$ is a continuous image of a supercompact space then the remainder $K {\setminus} \natN$ satisfies the countable chain condition. 
The Stone space of $B(\calA)$ is a compactification of $\omega$ 
whose the remainder is homeomorphic to $A(\kappa)$, and thus not c.c.c..
Therefore this space is not a continuous image of any supercompact space and 
in particular $B(\calA)$ is not embeddable into any poset algebra.
\begin{lemma}
\label{lemma-3.3}
Let $P$ be a poset such that $A(\aleph_1)$ embeds
topologically into $\Fs(P)$.
Then $P$ contains an uncountable antichain.
\end{lemma}
\begin{pf}
Assume $A(\aleph_1)$ is a closed subspace of $\Fs(P)$.
For each $\alpha \in \aleph_1$ pick a proper clopen prime ideal 
(or a clopen prime filter)
$H_\alpha$ of $\Fs(P)$ such that $\alpha \in H_{\alpha}$ and
$\infty \not\in H_\alpha$.
Replacing $\aleph_1$ by its uncountable subset, we may assume that all the sets $H_\alpha$ are of the same type, i.e. either all of them are prime ideals or all of them are prime filters.

For each $\alpha < \kappa$ the set $s_\alpha = H_\alpha \cap \aleph_1$ is finite, 
because $H_\alpha$ is a closed set which does not contain 
the accumulation point $\infty$ of $A(\aleph_1)$.
Furthermore $\alpha\in s_\alpha$.
By the $\Delta$-system Lemma,
there are an uncountable set $S\subseteq\aleph_1$ and a finite set
$t \subseteq \aleph_1$ such that $s_\alpha \cap s_\beta = t$
for distinct  $\alpha,\beta \in S$.
If $\alpha \neq \beta$ are in $S \setminus [0,\max(t))$ then
$H_\alpha \not\subseteq H_\beta$, because $\alpha\in H_\alpha {\setminus} H_\beta$.
It follows that $\setn{ H_\alpha }{ \alpha \in S \setminus [0,\max(t))}$
is an uncountable antichain with respect to inclusion.
By \cite[Lemma 2.9]{ABKR},
every proper clopen prime filter in $\Fs(P)$ is of the form
$V_p$ for some $p\in P$ and 
every proper clopen prime ideal is of the form $-V_p$ for some $p\in P$.
Since $p \leq q$ iff $V_p \subseteq V_q$,
in both cases $P$ contains an uncountable antichain.
\end{pf}
\begin{theorem}
\label{thm-3.4}
For every uncountable almost disjoint family 
$\calA \subseteq \calP(\omega)$ the algebra $B(\calA)$
is well generated and is not embeddable into any poset algebra.

Furthermore, if\/ $B(\calA)$ is a homomorphic image of\/ $\F(P)$, 
then $P$ has an uncountable antichain.
\end{theorem}
\begin{pf} To see that $B(\calA)$ is well generated,
note that the meet subsemilattice generated by 
$\calG =\setn{ \sngltn{n} }{ n\in\omega } \cup\calA$ is well founded
and thus, by Proposition~\ref{prop-2.3}(a2), the lattice generated by~$\calG$ is well-founded.
The fact that $B(\calA)$ is not embeddable into a poset Boolean
algebra follows from the result of Bell \cite[Corollary 3.2]{bell}
quoted above.
The ``furthermore" part follows from Lemma~\ref{lemma-3.3}:
if $\fnn{h}{\F(P)}{B(\calA)}$ is a homomorphism onto,
then the Stone space $Y$ of $B(\calA)$ is embeddable into $\Fs(P)$.
Removing the isolated points from $Y$ we get a copy of $A(\kappa)$,
where $\kappa = |\calA| > \aleph_0$.
Hence $A(\aleph_1)$ is topologically embeddable into $\Fs(P)$. 
By Lemma~\ref{lemma-3.3}, $P$ contains an infinite antichain.
\end{pf}

The next section shows that $\gbqo {\setminus} \WQO \neq \emptyset$.
\subsection{A poset Boolean algebra over a scattered chain,
generated by a better quasi-ordered set,
not isomorphic to any poset algebra over a well quasi-ordering}
\label{S3c}
The announced example is the algebra $\F(\omega_1+\omega^*)$.
We start with an auxiliary result on topological lattices.
\begin{lemma}
\label{lemma-3.5}
Let $\fifthtpl{L}{0^L}{1^L}{\meet}{\join}$ be a topological lattice,
and assume that:
\begin{enumerate}
\vspace{-2.5mm}
\item[\rfs{1.}] 
$1^L$ is a limit of a nontrivial convergent sequence;
\vspace{-2.5mm}
\item[\rfs{2.}]
there exists a sequence $\seqn{ a_\alpha }{\alpha < \omega_1 }$
such that $\lim_{\alpha<\omega_1} a_\alpha = 1^L$ and $a_\alpha<1^L$
for every $\alpha < \omega_1$.
\end{enumerate}
\vspace{-2.5mm}
Then $L$, as a topological space, is not linearly orderable.
\end{lemma}
\begin{pf} Let $\seqn{b_n}{n\in\omega}$ be a sequence converging to
$1^L$ such that $b_n < 1^L$ for every $n\in\omega$.
Fix $n\in\omega$. 
By~(b), we have
($\dagger$):~$\lim_{\alpha<\omega_1}(a_\alpha {\meet} b_n) 
= 1^L {\meet} b_n=b_n$.
Now suppose that $\prec$ is a strict linear order on $L$
inducing the topology of $L$.
Taking a subsequence if necessary and possibly reversing the order,
we may assume that $b_n \prec 1^L$ for every $n\in\omega$ and that
$b_n \prec b_m$ whenever $n<m$.
Using ($\dagger$), we can find $\alpha_n<\omega_1$ such that
($\ddagger$):~$b_{n-1} \prec a_\alpha {\meet} b_n \prec b_{n+1}$ for
$\alpha \geq \alpha_n$.
Let $\beta = \sup_{n\in\omega}\alpha_n$.
We have
$\lim_{n} \,a_\beta {\meet} b_n =a_\beta {\meet} 1^L = a_\beta$.
However, by ($\ddagger$),
$b_{n-1} \prec a_\beta {\meet} b_n \prec b_{n+1}$ for every $n$,
and thus $\lim_{n} \, a_\beta {\meet} b_n = 1^L$.
This is a contradiction.
\end{pf}
\begin{lemma}\label{lemma-3.6}
Let $P$ be a well quasi-ordered poset.
Then every $z \in \Fs(P)$ is isolated in 
$\setn{ x \in \Fs(P) }{ x \subseteq z }$.
\end{lemma}
\begin{pf}
Let $\sigma$ consist of all minimal elements of $z$.
Then $\sigma$ is finite and
$V \eqdef \setn{ y \in \Fs(P) }{ \sigma \subseteq y }$
is a neighborhood of $z$ which is disjoint from 
$\setn{ y \in \Fs(P) }{ y \subseteq z \mbox{\rm\ and } y \neq z}$.
\end{pf}
\begin{theorem}
\label{thm-3.7}
$\F(\omega_1+\omega^*)$ is generated by a better quasi-ordered
subset, and this algebra is not isomorphic to $\F(P)$ 
for any well quasi-ordered poset $P$.
\end{theorem}
\begin{pf}
To see that $\F(\omega_1+\omega^*)$ is bqo generated,
note that
$G=\setn{ x_{\alpha} }{ \alpha<\omega_1 } \cup
\setn{ -x_n }{ n\in\omega^* }$
is isomorphic to the disjoint sum of $\omega_1$ and $\omega$
and therefore it is better quasi-ordered.
Clearly, $G$ generates $\F(\omega_1+\omega^*)$.

The Stone space of $\F(\omega_1+\omega^*)$, which we shall denote
by $X$, is homeomorphic to the linearly ordered space
$\omega_1+1+\omega^*$. 
Suppose $\F(\omega_1+\omega^*)$ is isomorphic to $\F(P)$ for 
some wqo set $P$.
Topologically, this means that $X$ is isomorphic to a well founded,
compact, 0-dimensional, distributive topological lattice 
$\fifthtpl{L}{0^L}{1^L}{\meet}{\join}$, where $L = \Fs(P)$
and the order of $L$ is the reversed inclusion.
Furthermore $L = A \cup \sngltn{p} \cup B$,
where $A \cup \sngltn{p}$ is homeomorphic to $\omega_1+1$
and $B \cup \sngltn{p}$ is homeomorphic to $\omega+1$
(and $p$ is the limit point in both cases). 

Let $\preceq$ denote the lattice ordering of $L$.
By Lemma \ref{lemma-3.6}, $p$ is isolated in $[p,1^L]^{\preceq}$,
which implies that $a {\meet} p \prec p$ for all 
but countably many $a \in A$ and $b {\meet} p \prec p$ for all
but finitely many $b \in B$.
In particular, 
$p$ is the limit of a non-trivial sequence in $[0^L,p]^{\preceq}$.
Enumerating suitably the set
$\setn{ a {\meet} p }{ a\in A \mbox{\rm\ and } a {\meet} p \prec p}$,
we see that the assumptions of Lemma \ref{lemma-3.5} are satisfied
for the lattice $[0^L,p]^{\preceq}$.
By this Lemma, $[0^L,p]^{\preceq}$ (and therefore also $L$) 
is not linearly orderable, which is a contradiction.
\end{pf}
\subsection{A Boolean algebra generated by a better quasi-ordered
subset and not isomorphic to any poset algebra}
\label{S3d}
The announced example is the clopen algebra of the compact space $X$
obtained from three disjoint copies of $\omega_1+1$ by identifying
the three complete accumulation points.
Note that $X$ is homeomorphic to the subspace of $(\omega_1+1)^2$
which is the union of the diagonal, horizontal and vertical lines
passing through the point $\pair{\omega_1}{\omega_1}$.
\begin{lemma}
\label{lemma-3.8}
Let $\fifthtpl{L}{0^L}{1^L}{\meet}{\join}$ be a topological lattice
and assume that there exist two sequences 
$\seqn{a_\alpha}{\alpha<\omega_1}$
and $\seqn{b_\beta}{\beta<\omega_1}$ such that:
\begin{enumerate}
\vspace{-2mm}
\item[\rfs{(1)}] $\seqn{a_\alpha}{\alpha<\omega_1}$ and
$\seqn{b_\beta}{\beta<\omega_1}$
are continuous, i.e. $a_\lambda = \lim_{\alpha<\lambda} a_\alpha$
and $b_\lambda = \lim_{\alpha<\lambda} b_\alpha$ for every limit
$\lambda<\omega_1$.
\vspace{-2mm}
\item[\rfs{(2)}] $\lim_{\alpha<\omega_1} a_\alpha = 1^L 
= \lim_{\alpha<\omega_1} b_\alpha$
and $a_\alpha, b_\alpha <1^L$ for every $\alpha<\omega_1$.
\vspace{-2mm}
\item[\rfs{(3)}] For every $\alpha,\beta<\omega_1$, $a_\alpha$ and
$b_\beta$ are $G_\delta$-points in $L$.
\vspace{-1mm}
\end{enumerate}
Then there exists a closed unbounded set $C\subseteq\omega_1$ such that $a_\alpha = b_\alpha$ for $\alpha \in C$.
\end{lemma}
\begin{pf} 
Fix $\beta<\omega_1$.
By~(2) we have 
($*$):~$\lim_{\alpha<\omega_1}(a_\alpha {\meet} b_\beta) 
= 1^L {\meet} b_\beta = b_\beta$.
By~(3) and ($*$), for each $\beta<\omega_1$ there exists $\alpha(\beta)<\omega_1$ such that
$a_\alpha {\meet} b_\beta = b_\beta$,
that is $b_\beta \leq a_\alpha$, for $\alpha \geq \alpha(\beta)$.
By symmetry, for each $\alpha<\omega_1$ there exists
$\beta(\alpha)<\omega_1$ such that $a_\alpha \leq b_\beta$
whenever $\beta \geq \beta(\alpha)$.
Let $C$ consist of all limit ordinals $\lambda < \omega_1$ such that
$\alpha(\xi) < \lambda$ and $\beta(\xi) < \lambda$
for every $\xi < \lambda$.
By~(1), $a_\lambda = b_\lambda$ for every $\lambda \in C$.
\end{pf}
\begin{lemma}
\label{lemma-3.9}
Let $X$ be the space obtained by taking three disjoint copies of
$\omega_1+1$ and identifying the last points.
Then $X$ does not have a structure of a topological lattice.
\end{lemma}
\begin{pf}
Write $X=\sngltn{p} \cup A\cup B\cup C$, where $A=\setn{a_\alpha}{\alpha<\omega_1}$, $B=\setn{b_\beta}{\beta<\omega_1}$ and $C=\setn{c_\alpha}{\alpha<\omega_1}$, the enumerations are
continuous, and 
$p = \lim_{\alpha<\omega_1} a_\alpha 
= \lim_{\alpha<\omega_1} b_\alpha 
= \lim_{\alpha<\omega_1} c_\alpha$.
Note that $X$ is $0$-dimensional and compact, and
that $p$ is the only point of $X$ with uncountable character. 

Suppose $\fifthtpl{L}{0^L}{1^L}{\meet}{\join}$ is a topological lattice and denote by $\preceq$ the lattice order of $X$.
It cannot be the case that both $A \cap [p,1^X]^{\preceq}$ 
and $B \cap [p,1^X]^{\preceq}$
are uncountable, because then by Lemma~\ref{lemma-3.8}, 
we would get a contradiction.
So assume that $A \cap [p,1_X]^{\preceq}$ is countable and also that
$C \cap [p,1^X]^{\preceq}$ is countable. 
Without loss of generality, we can assume that
$[p,1^X]^{\preceq} \subseteq \sngltn{p} \cup B$, 
by adding to $B$ a closed countable set.
Similarly, 
assume that $[0^X,p]^{\preceq} \subseteq \sngltn{p} \cup A$.
Now, every element of $C$, except possibly countably many,
must be incomparable with $p$, since otherwise we could use 
Lemma~\ref{lemma-3.8} again to get a contradiction.
So assume $c_\alpha$ is not comparable with $p$ for every
$\alpha<\omega_1$. 

We have 
$\lim_{\alpha<\omega_1}\,c_\xi {\meet} c_\alpha 
= c_\xi{\meet} p\in A$.
Thus, by the fact that each point of $A$ is $G_\delta$,
for every $\xi < \omega_1$ we can find $\alpha(\xi) > \xi$ such that
$c_\xi {\meet} c_\alpha \in A$ for $\alpha \geq \alpha(\xi)$
(in fact $c_\xi {\meet} c_\alpha$ is eventually constant 
with respect to $\alpha$).

Finally, choose an increasing sequence $\seqn{\xi_n}{n\in\omega}$
such that $\alpha(\xi_n)<\xi_{n+1}$ for every $n\in\omega$.
Let $\beta= \lim_{n}\xi_n$.
Then
$c_\beta = c_\beta {\meet} c_\beta 
= \lim_{n}(c_{\xi_n} {\meet} c_{\alpha(\xi_n)})$, and
so\break
$c_\beta\in A$, because $A$ is closed under countable limits.
This is a contradiction.
\end{pf}
\begin{theorem}
\label{thm-3.10}
There exists a Boolean algebra $B$ generated by a 
better quasi-ordered subset,
namely a homomorphic image of the poset algebra of 
the disjoint sum of two copies of\/  $\omega_1$, 
which is not isomorphic to a poset algebra.
\end{theorem}
\begin{pf}
We use duality. 
The space $X$ from Lemma~\ref{lemma-3.9} is homeomorphic
to a closed subspace of $(\omega_1+1){\times}(\omega_1+1)$ 
and the latter space corresponds to the poset algebra of 
the disjoint sum of two copies of $\omega_1$.
\end{pf}

\end{document}